\documentclass[10pt,a4paper]{article}

\usepackage[dvips]{color}

\usepackage{times}

\usepackage{macros-pdcat08}

\title{A Data-Parallel Algorithm to Reliably Solve\\Systems of Nonlinear Equations}
 \author{Fr\'ed\'eric Goualard${}^{\P,}\footnotemark[6]\;$ and Alexandre Goldsztejn${}^{\S}$\\[5pt]
\normalsize${}^\P$Universit\'e de Nantes, Nantes Atlantique Universit\'es, CNRS, LINA, UMR 6241.\\
\normalsize2 rue de la Houssini\`ere, BP 92208, F-44000 NANTES \\
\normalsize${}^\S$CNRS, LINA, UMR 6241. 2 rue de la Houssini\`ere, BP 92208, F-44000 NANTES}
\date{}

\begin{document}

\maketitle
\footnotetext[6]{Corresponding author's email: \texttt{Frederic.Goualard@univ-nantes.fr}.}

\begin{abstract}
\noindent Numerical methods based on interval arithmetic are efficient means to reliably 
solve nonlinear systems of equations. Algorithm bc3revise is an interval method that 
tightens variables' domains by enforcing a property called 
box consistency. It has been successfully used on difficult problems whose solving 
eluded traditional numerical methods.
We present a new algorithm to enforce box consistency that is simpler than bc3revise, 
faster, and easily data parallelizable. A parallel implementation with Intel SSE2 SIMD instructions 
shows that an increase in performance of up to an order of magnitude and more is achievable.

\end{abstract}

\section{Introduction}

Interval methods~\cite{Neumaier:Book1990} are numerical algorithms that use \emph{interval arithmetic}~\cite{Moore:Book1966}
to avoid rounding error problems intrinsic to floating-point arithmetic~\cite{Higham:2002:ASN}. 
They give enclosures of all solutions of nonlinear systems of equations with the guarantee that no solution is 
ever lost.

Straight interval extensions of classical numerical algorithms such as the Newton method are not well-suited to 
problems with many solutions or with large initial domains for the variables. To tackle these shortcomings, 
elaborate algorithms have been devised in the context of \emph{Interval Constraint Programming}~\cite{Benhamou:EO2001};
they are usually employed as the inner stage of a free-steering nonlinear Gauss-seidel method 
to exclude parts of a variable's domain that do not contain zeroes of a unidimensional equation. Domain tightening
is achieved by enforcing some local consistency property. \emph{Box consistency}~\cite{Benhamou-et-al:ILPS94} is one 
such consistency notion, which has been proved efficient in handling hard problems whose solving eluded traditional numerical 
methods for years~\cite{Granvilliers_Benhamou:JGO2001}. It is usually enforced by Algorithm
 \FunT{bc3revise}~\cite{Benhamou-et-al:ILPS94}, which 
combines a binary search with interval Newton steps~\cite{Moore:Book1966} to isolate leftmost and rightmost zeroes of a 
unidimensional equation in the domain of a variable. 

Thanks to ubiquitous Intel SSE2 SIMD instructions, it is possible to perform many interval 
operations at roughly the same cost as floating-point operations by computing the
two bounds of the result in parallel (\emph{basic interval vectorization})~\cite{Goualard:PARA2008}. 
We outline in Section~\ref{sec:interval-arithmetic} a novel way to do even better and to compute an interval 
function for two different intervals in parallel (a four times speed-up compared to ``sequential'' interval
evaluation).

As all interval methods, Algorithm \FunT{bc3revise}---described in Section~\ref{sec:box-consistency-and-bc3}---can 
benefit from basic interval vectorization without any
modification. On the other hand, it cannot take full advantage of the new arithmetic described in Section~\ref{sec:interval-arithmetic}. 
Hence the introduction of Algorithm \FunT{sbc} in Section~\ref{sec:box-consistency-by-shaving}: 
it is a new algorithm that enforces box consistency by ``shaving'' domains from the left and right bounds 
inward. Experiments show that a sequential version of \FunT{sbc} is already faster than \FunT{bc3revise}
on a set of test problems. We then describe in Section~\ref{sec:simd-algorithm} an algorithm  
that exploits the potential for a high level of data parallelism in \FunT{sbc} by using SSE2 instructions
to perform interval arithmetic evaluations of functions at four times the speed of a sequential code.
Experiments are reported in Section~\ref{sec:experiments} and show increases in performances  over 
\FunT{bc3revise} of up to an order of magnitude and more.

\section{Interval Arithmetic and its Vectorization}
\label{sec:interval-arithmetic}

Classical iterative numerical methods suffer from defects such as loss
of solutions, absence of convergence, and convergence to unwanted
attractors due to the use of floating-point numbers
(aka \emph{floats}). At the end of the fifties, Moore~\cite{Moore:Book1966}
popularized the use of intervals to control the errors made while
computing with floats. Additionnally, interval extensions of iterative numerical methods are
always convergent.

In the following, we use the notations sponsored by Kearfott and others~\cite{Kearfott-et-al:Baikal2005}, where interval
quantities are boldfaced.

Interval arithmetic replaces floating-point numbers
by closed connected sets of the form
$\itvK{I}=\itv{\lbK{\itvK{I}}}{\rbK{\itvK{I}}}=\{a\in\RSet\mid
\lbK{\itvK{I}}\leq a\leq\rbK{\itvK{I}}\}$ from \emph{the set \ISet\ of intervals},
where \lbK{\itvK{I}} and
\rbK{\itvK{I}} are floating-point numbers. In addition, each $n$-ary real
function $\phi$ with domain $\Domain{\phi}\subseteq\RSet^n$ is \emph{extended to an interval
function} $\Phi$ with domain $\Domain{\Phi}\subseteq\ISet^n$ in such a way
that the \emph{containment principle} is verified:
\begin{equation}\label{eq:containment-principle}
  \forall \vecK{a}\in\Domain{\phi}, \forall\vecK{\itvK{I}}\in\Domain{\Phi}\colon
  \vecK{a}\in\vecK{\itvK{I}}\implies\mathop{\phi}(\vecK{a})\in\mathop{\Phi}(\vecK{\itvK{I}}),
\end{equation}
as illustrated by the following example.

\begin{example}
  The \emph{natural interval extensions} of addition and
  multiplication are defined by:
\begin{align*}
  \itvK{I_1}+\itvK{I_2}&=\itv{\RndDn{\lbK{\itvK{I_1}}+\lbK{\itvK{I_2}}}}{\RndUp{\rbK{\itvK{I_1}}+\rbK{\itvK{I_2}}}}\\
  \itvK{I_1}\times\itvK{I_2}&=
 \itv{\min(\RndDn{\lbK{\itvK{I_1}}\lbK{\itvK{I_2}}},\RndDn{\lbK{\itvK{I_1}}\rbK{\itvK{I_2}}},\RndDn{\rbK{\itvK{I_1}}\lbK{\itvK{I_2}}},\RndDn{\rbK{\itvK{I_1}}\rbK{\itvK{I_2}}})}{\max(\RndUp{\lbK{\itvK{I_1}}\lbK{\itvK{I_2}}},\RndUp{\lbK{\itvK{I_1}}\rbK{\itvK{I_2}}},\RndUp{\rbK{\itvK{I_1}}\lbK{\itvK{I_2}}},\RndUp{\rbK{\itvK{I_1}}\rbK{\itvK{I_2}}})}
\end{align*}
where \RndDn{r} (resp., \RndUp{r}) is the greatest floating-point number
smaller or equal (resp., the smallest floating-point number greater or equal) to $r$. 

Then, given the real function $f(x,y)=x\times x+y$, we may define its natural interval extension by
$\itvK{f}(\itvK{x},\itvK{y})=\itvK{x}\times\itvK{x}+\itvK{y}$, and we have that, e.g.,
$\itvK{f}(\itv{2}{3},\itv{-1}{5})=\itv{3}{14}$.
\end{example}
Implementations of interval arithmetic use outward rounding to enlarge
the domains computed so as not to violate the containment principle~\eqref{eq:containment-principle},
should some bounds be unrepresentable
with floating-point numbers.

Interval addition, subtraction, multiplication, division, and integral exponentiation may be computed at
roughly the same speed as their floating-point counterpart thanks to SIMD instructions, and in particular, to
Intel SSE2 instructions.

\begin{figure}[htbp]
\begin{center}
\subfigure[SIMD floating-point arithmetic]{\label{fig:simd-arithmetic}%
	\includegraphics[scale=.4]{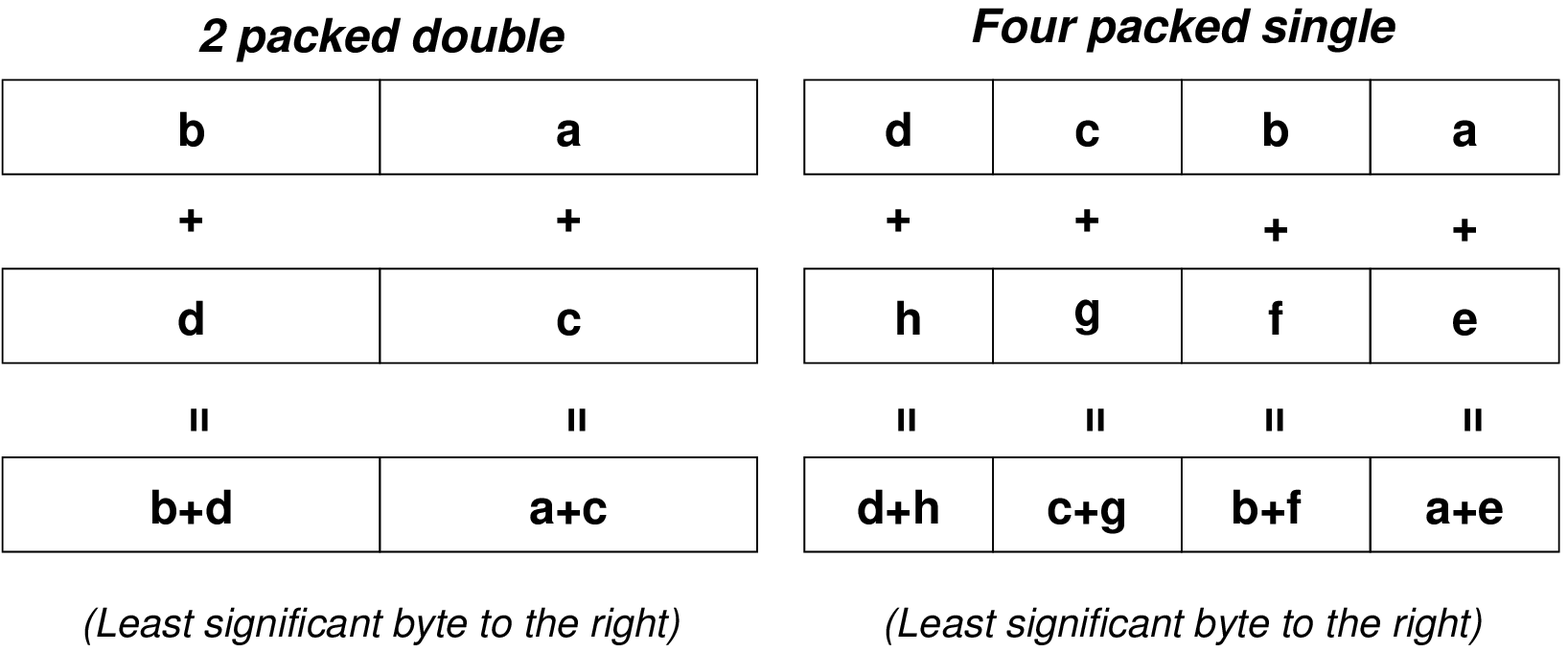}}\hspace{20pt}
\subfigure[Two interval additions with one SSE2 instruction]{\label{fig:simd-interval-addition}%
  	\includegraphics[scale=.4]{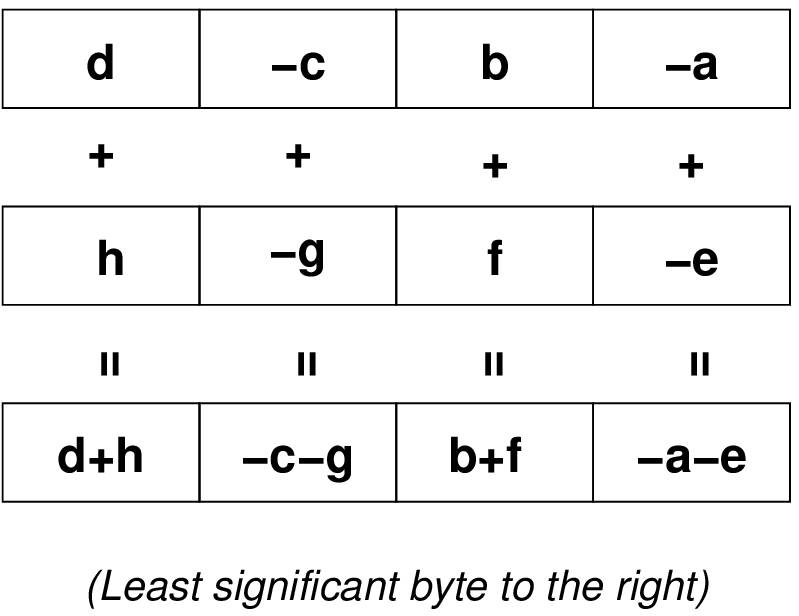}
}
\caption{Floating-point arithmetic and interval arithmetic in SSE2 registers}
\label{fig:simd-registers}
\end{center}
\end{figure}

Intel SSE2 instructions manipulate 128 bits registers that may be interpreted in various ways. Most notably, the registers
may pack 2 double precision or 4 single precision floating-point numbers. An SSE2 operator may then compute 2 or 4 
floating-point operations in parallel (see Figure~\ref{fig:simd-arithmetic}).

The direction of rounding for SSE2 instructions is selected independently of that of the Floating-Point Unit (FPU). 
An SSE2 instruction uses the same rounding for all operations
performed in parallel. Nevertheless, thanks to simple floating-point properties, it is still possible to write algorithms that compute 
in parallel the two outward-rounded bounds of the result of interval operations. For example, we may use the property:
\begin{equation*}
 \RndDn{a+b} = -\RndUp{-a-b}
\end{equation*}
where $a$ and $b$ are floating-point numbers. 

By \emph{storing the negation of the left bound} of an interval, and by setting once and for all the rounding direction for SSE2 instructions 
to $+\infty$, the two interval additions $\itv{a}{b}+\itv{e}{f}$ and $\itv{c}{d}+\itv{g}{h}$ can be performed by the sole SIMD 
instruction depicted in Figure~\ref{fig:simd-interval-addition}.
All the other operators may be defined accordingly. Goualard's paper~\cite{Goualard:PARA2008} illustrates these principles for the
case of basic interval vectorization (two double precision bounds computed in parallel). The algorithms to compute four bounds
in parallel are new and are reported in an unpublished paper currently under review.

Armed with an interval library whose operators compute two interval operations in parallel, we may evaluate the interval extension of a function \itvK{f}
for two different intervals for the same cost as \emph{one floating-point} evaluation of $f$. In the following, we note \simd{\itvK{f}(\itvK{I_1})}{\itvK{f}(\itvK{I_2})}
such a parallel evaluation of \itvK{f} for two different interval arguments.

\section{Box Consistency and the bc3revise Algorithm}
\label{sec:box-consistency-and-bc3}

\emph{Interval Constraint Programming}~\cite{Benhamou:EO2001}
is a successful approach to reliably isolate all solutions of 
systems of equations. It makes cooperate
\emph{contracting operators} to prune the domains of the variables (intervals with floating-point bounds from the set \ISet) with
smart \emph{propagation algorithms}~\cite{Mackworth:AI77}---akin to free-steering nonlinear Gauss-Seidel---to ensure consistency among all the constraints.

The amount of pruning obtained from one equation is controlled by the level of consistency enforced. 
Box consistency~\cite{Benhamou-et-al:ILPS94} is defined as follows:
\begin{definition}[Box consistency]
 An equation $f(x_1,\dots,x_n)=0$ is box consistent with respect to a variable $x_i$ and 
 a box $\vecK{\itvK{B}}=\itvK{I_1}\times\cdots\times \itvK{I_n}$ if and only if:
 \begin{equation}\label{eq:box-consistency}
  \left\{\begin{array}{l}
	0\in\itvK{f}(\itvK{I_1},\dots,\itvK{I_{i-1}},\itv{\lbK{\itvK{I_i}}}{\nextFloat{\lbK{\itvK{I_i}}}},\itvK{I_{i+1}},\dots,\itvK{I_n})       \\
	\text{and}\\
  	0\in\itvK{f}(\itvK{I_1},\dots,\itvK{I_{i-1}},\itv{\previousFloat{\rbK{\itvK{I_i}}\,}}{\rbK{\itvK{I_i}}},\itvK{I_{i+1}},\dots,\itvK{I_n}), \\
  \end{array}\right.
 \end{equation}
where $\itvK{I}=\itv{\lbK{\itvK{I}}}{\rbK{\itvK{I}}}$ is an interval with floating-point bounds, \nextFloat{a} (resp., \previousFloat{a}) is
the smallest floating-point number greater than (resp., the largest floating-point number smaller than) the floating-point number $a$, 
and \itvK{f} is the natural interval extension of $f$.
\end{definition}

Given a real function $f\colon\RSet^n\to\RSet$, and a box of domains 
$\itvK{B}=\itvK{I_1}\times\cdots\times \itvK{I_n}\in\ISet^n$, we define $\itvK{g^B_i}:\ISet\to\ISet$ as 
the $i$th unary interval projection with respect to \itvK{B} of its interval extension \itvK{f}:
\begin{equation*}
 \itvK{g^B_i}(x) = \itvK{f}(\itvK{I_1},\dots,\itvK{I_{i-1}},x,\itvK{I_{i+1}},\dots,\itvK{I_n}), \quad i\in\{1,\dots n\}.
\end{equation*}
In the following, we will mostly manipulate $\itvK{g^B_i}$ instead of $\itvK{f}$. The original real function $f$, 
the box $\itvK{B}$ of domains considered and/or the variable $x_i$ on which the projection is performed will often be left implicit
and omitted from the notation of \itvK{g}. 

In order not to lose any potential solution, an algorithm that enforces box consistency of an equation
with respect to a variable $x_i$ and a box of domains must return the largest domain $\itvK{I'_i}\subseteq \itvK{I_i}$
that verifies Eq.~\eqref{eq:box-consistency}.

\begin{algorithm}[htbp]
 \begin{center}
  \caption{Computing a box consistent interval with respect to \itvK{g} the usual way}
  \label{alg:bc3}
  \algoInitCounter
  \begin{tabbing}
   \hspace*{1cm}\=12\=12\=12\=\kill
   \textbf{[\FunT{bc3revise}]} \KWD{in:} $\itvK{g}\colon\ISet\to\ISet$; \KWD{in:} $\itvK{I}\in\ISet$\\
   \REM{Returns the largest interval included in \itvK{I} that is box consistent with respect to \itvK{g}}\\
   \KWD{begin}\\
   \ACTR\> $\itvK{I_l}\gets\Fun{left\_narrow}(\itvK{g},\itvK{I})$\\
   \ACTR\>\KWD{if} $\itvK{I_l} \neq\emptyset\colon$\\
   \ACTR\>\>$\itvK{I_r}\gets\Fun{right\_narrow}(\itvK{g},\itv{\lbK{I_l}}{\rbK{I}})$\\
   \ACTR\>\>\KWD{return} $\Hull(\itvK{I_l}\cup\itvK{I_r})$ \REM{returns the smallest interval w.r.t. 
										set inclusion that contains $\itvK{I_l}\cup\itvK{I_r}$}\\
   \ACTR\>$\KWD{else}\colon$\\
   \ACTR\>\>\KWD{return} $\emptyset$\\
   \KWD{end}
  \end{tabbing}
 \end{center}
\end{algorithm}

Algorithm~\FunT{bc3revise}~\cite{Benhamou-et-al:ILPS94}, presented by Algorithms~\ref{alg:bc3} 
and ~\ref{alg:left_narrow}, considers the unary projection of an $n$-ary equation on a variable $x_i$
(where all variables but $x_i$ have been replaced by their current domain) and a domain $\itvK{I_i}$.
It enforces box consistency by searching the leftmost and rightmost \emph{canonical domains}\footnote{A non-empty interval 
\itv{a}{b} is canonical if $\nextFloat{a}\geq b$.} for which \itvK{g} evaluates to an interval containing $0$.
The search is performed by a dichotomic search aided with Newton steps to accelerate the process. Algorithm~\ref{alg:left_narrow}
describes the search of a quasi-zero to update the left side of $\itvK{I_i}$. The procedure \FunT{right\_narrow} to update the right bound
works along the same lines and is, therefore, omitted.

\begin{algorithm}[htbp]
 \begin{center}
  \caption{Computing a box consistent left bound with Newton steps and a binary search}
 \label{alg:left_narrow}
 \algoInitCounter
 \begin{tabbing}
  \hspace*{1cm}\=12\=12\=12\=12\=\kill
  \textbf{[\FunT{left\_narrow}]} \KWD{in:} $\itvK{g}\colon\ISet\to\ISet$; \KWD{in:} $\itvK{I}\in\ISet$\\
   \REM{Returns an interval included in \itvK{I}}\\
   \REM{with the smallest left bound $l$ such that $0\in\itvK{g}(\itv{l}{\nextFloat{l}})$}\\
 \KWD{begin}\\
 \ACTR\>\KWD{if} $0\notin\itvK{g}(\itvK{I})\colon$ \REM{No solution in \itvK{I}}\\
 \ACTR\>\> return $\emptyset$\\
 \ACTR\>\KWD{else}:\\
 \ACTR\>\>\KWD{if} $\nextFloat{\lbK{I}}\geq\rbK{I}$: \REM{\Fun{canonical}(\itvK{I})}:\\
 \ACTR\>\>\>\KWD{return} \itvK{I}\\
 \ACTR\>\>\KWD{else}:\\
 \ACTR\>\>\>$\itvK{I}\gets\Fun{Newton}(\itvK{g},\itvK{g'},\itvK{I})$ \REM{Interval Newton steps}\\
 \ACTR\>\>\>\KWD{if} $0\in\itvK{g}(\itv{\lbK{I}}{\nextFloat{\lbK{I}}})$: \REM{Box consistent left bound?}\\
 \ACTR\>\>\>\>\KWD{return} $\itvK{I}$\\
 \ACTR\>\>\>\KWD{else}:\\
 \ACTR\>\>\>\>$(\itvK{I_1},\itvK{I_2})\gets\Fun{split}(\itvK{I})$ \REM{$\itvK{I_1}\gets\itv{\lbK{I}}{\midpointK{I}}, \itvK{I_2}\gets\itv{\midpointK{I}}{\rbK{I}}$}\\
 \ACTR\>\>\>\>$\itvK{I}\gets\Fun{left\_narrow}(\itvK{g},\itvK{I_1})$\\
 \ACTR\>\>\>\>\KWD{if} $\itvK{I}=\emptyset\colon$\\
 \ACTR\>\>\>\>\>\KWD{return} $\Fun{left\_narrow}(\itvK{g},\itvK{I_2})$\\
 \ACTR\>\>\>\>\KWD{else}:\\
 \ACTR\>\>\>\>\>\KWD{return} $\itvK{I}$\\
 \KWD{end}
 \end{tabbing}
 \end{center}
\end{algorithm}

\begin{figure}[htbp]
 \begin{center}
	\includegraphics[bb=56 486 332 768]{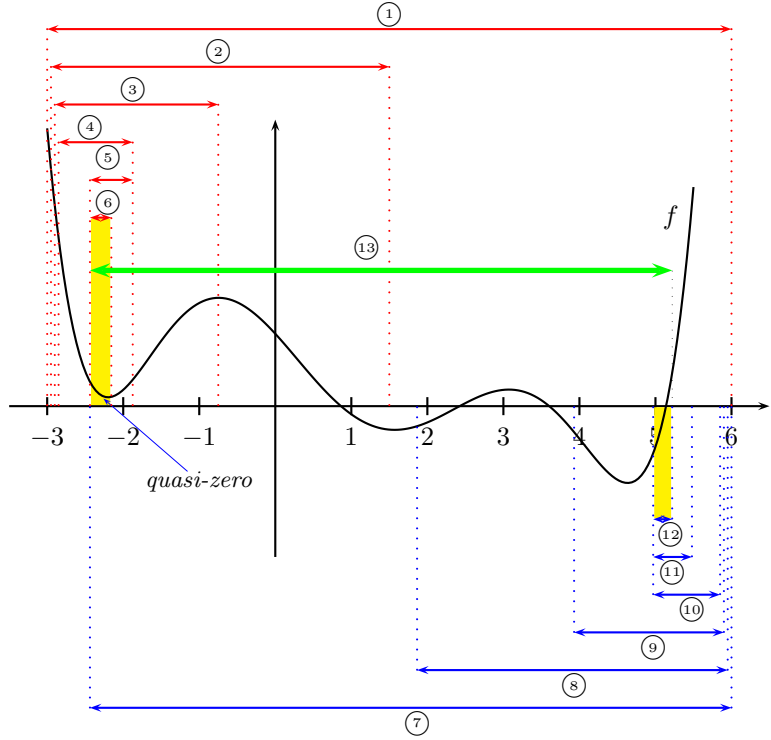}
  \caption{Enforcing box consistency with \FunT{bc3revise}}
	\label{fig:bc3-computation}
 \end{center}
\end{figure}

Algorithm~\FunT{bc3revise} first tries to move the left bound of $\itvK{I_i}$ to the right, and then proceeds to move its right bound 
to the left. The \FunT{Newton} procedure computes a fixpoint of the Interval Newton algorithm~\cite{Moore:Book1966}, where a Newton step at iteration $j+1$ is:
\begin{equation*}
 \itvK{I^{(j+1)}} \gets \itvK{I^{(j)}}\cap\left(\midpointK{I^{(j)}} - \frac{\itvK{g}(\midpointK{I^{(j)}})}{\itvK{g'}(\itvK{I^{(j)}})}\right)
\end{equation*}
with \midpointK{I} the midpoint of the interval \itvK{I}\footnote{Note that we are free to
choose any point in \itvK{I}, not the midpoint only. We take advantage of this in 
the algorithms presented in the next section.}. As in Ratz's work~\cite{Ratz:TR1996}, the Newton step uses an extended version
of the interval division to return a union of two semi-open-ended intervals whenever $\itvK{g'}(\itvK{I^{(j)}})$ contains $0$. The subtraction and
the intersection operators are modified accordingly. The intersection operator is applied to an interval (\itvK{I^{(j)}}) and 
a union of two intervals (result of the subtraction), and returns an interval.
Figure~\ref{fig:bc3-computation} presents graphically the steps performed to enforce box consistency. The encircled numbers
label the steps.

Algorithms \FunT{bc3revise}, \FunT{left\_narrow} and \FunT{right\_narrow} do not offer opportunities to exploit full 
data parallelism as they do not require close
evaluations of the same function over different domains. The same holds for the \FunT{Newton} procedure: in the general case, \itvK{g}
and \itvK{g'} are different functions, and therefore cannot be evaluated in parallel with SIMD instructions.

\section{Box Consistency by Shaving}

To obtain a higher level of data parallelism, we propose a new algorithm to enforce box consistency on
the projection $\itvK{g^B_i}$ of an equation $f=0$ on a variable $x_i$: starting from the original domain $\itvK{I_i}$ for $x_i$,
consider separately its left half and its right half; for the left part, linearize  \itvK{g} at $\lbK{\itvK{I_i}}$ as
for a Newton step, solve the resulting linear equation and intersect the resulting domain with the left half of $\itvK{I_i}$; do the same
for the right half by linearizing \itvK{g} at $\rbK{\itvK{I_i}}$. A new smaller domain that preserves solutions is then obtained; Iterate
until reaching a fixpoint.

\subsection{A Sequential Algorithm: \FunT{sbc}}
\label{sec:box-consistency-by-shaving}

\begin{figure}[htbp]
\begin{center}
\includegraphics[scale=.8]{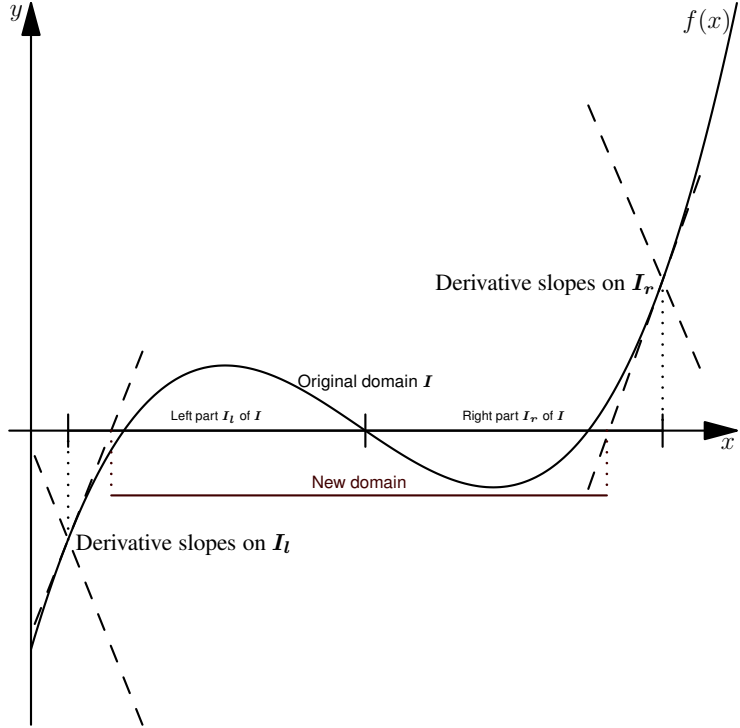}
 \caption{Domain reduction obtained with one iteration of the loop in \FunT{sbc}}
 \label{fig:sbc}
 \end{center}
\end{figure}

Figure~\ref{fig:sbc} illustrates graphically the process just described, and Algorithm~\ref{alg:sbc} presents the actual algorithm.

\begin{algorithm}[h!]
 \begin{center}
  \caption{Enforcing box consistency by shaving}
  \label{alg:sbc}
 \algoInitCounter
 \begin{tabbing}
  \hspace*{1cm}\=12\=12\=12\=12\=12\=\kill
  \textbf{[\FunT{sbc}]} \KWD{in:} $\itvK{g}\colon\ISet\to\ISet$; \KWD{in:} $\itvK{I}\in\ISet$\\[3pt]
   \REM{Returns the largest interval included in \itvK{I} that is box consistent w.r.t.\ \itvK{g}}\\[3pt]
 \KWD{begin}\\[3pt]
 \ACTR\>$(\Var{left\_consistent},\Var{right\_consistent})\gets(\KWD{false},\KWD{false})$\\[3pt]
 \ACTR\>\KWD{do}:\\[3pt]
 \ACTR\>\>$(\itvK{I_l},\itvK{I_r})\gets\Fun{split}(\itvK{I})$\\[3pt]
 \ACTR\>\>\KWD{if} $\neg\Var{left\_consistent}$: \REM{Updating the left bound}\\[3pt]
 \ACTR\>\>\>\KWD{if} $0\notin\itvK{g}(\itv{\lbK{I_l}}{\nextFloat{\lbK{I_l}}})$: \REM{\itvK{I} not box consistent to the left?}\\[3pt]
 \ACTR\>\>\>\>$\itvK{I_l}\gets\itv{\nextFloat{\lbK{I_l}}}{\rbK{I_l}}$ \REM{Considering the remainder of \itvK{I_l}}\\[3pt]
 \ACTR\>\>\>\>\KWD{if} $0\notin\itvK{g}(\itvK{I_l})$: \REM{No solution in the remainder of \itvK{I_l}?}\\[3pt]
 \ACTR\>\>\>\>\>$\itvK{I_l}\gets\emptyset$\\[3pt]
 \ACTR\>\>\>\>\KWD{else}:\\[3pt]
 \ACTR\>\>\>\>\>$\itvK{I_l}\gets \itvK{I_l}\cap\left(\lbK{I_l}-\itvK{g}(\itv{\lbK{I_l}}{\lbK{I_l}})/\itvK{g'}(\itvK{I_l})\right)$ \REM{One Newton step}\\[3pt]
 \ACTR\>\>\>\KWD{else}:\\
 \ACTR\>\>\>\>$\Var{left\_consistent}\gets\KWD{true}$\\
 \ACTR\>\>\KWD{if} $\neg\Var{right\_consistent}$: \REM{Updating the right bound}\\[3pt]
 \ACTR\>\>\>\KWD{if} $0\notin\itvK{g}(\itv{\previousFloat{\rbK{I_r}}}{\rbK{I_r}})$: \REM{\itvK{I} not box consistent to the right?}\\[3pt]
 \ACTR\>\>\>\>$\itvK{I_r}\gets\itv{\lbK{I_r}}{\previousFloat{\rbK{I_r}}}$ \REM{Considering the remainder of \itvK{I_r}}\\[3pt]
 \ACTR\>\>\>\>\KWD{if} $0\notin\itvK{g}(\itvK{I_r})$: \REM{No solution in the remainder of \itvK{I_r}?}\\[3pt]
 \ACTR\>\>\>\>\>$\itvK{I_r}\gets\emptyset$\\[3pt]
 \ACTR\>\>\>\>\KWD{else}:\\[3pt]
 \ACTR\>\>\>\>\>$\itvK{I_r}\gets \itvK{I_r}\cap\left(\rbK{I_r}-\itvK{g}(\itv{\rbK{I_r}}{\rbK{I_r}})/\itvK{g'}(\itvK{I_r})\right)$ \REM{One Newton step}\\[3pt]
 \ACTR\>\>\>\KWD{else}:\\
 \ACTR\>\>\>\>$\Var{right\_consistent}\gets\KWD{true}$\\
 \ACTR\>\>$\itvK{I}\gets\Hull(\itvK{I_l}\cup\itvK{I_r})$ \REM{Returns the ``hull'' $\itv{\min(\lbK{I_l},\lbK{I_r})}{\max(\rbK{I_l},\rbK{I_r})}$ of the union}\\
 \ACTR\>\KWD{while} $\left(\itvK{I}\neq\emptyset\right)\wedge\left(\neg\Var{left\_consistent}\vee\neg\Var{right\_consistent}\right)$\\[3pt]
 \ACTR\>\KWD{return} \itvK{I}\\[3pt]
  \KWD{end}
  \end{tabbing}
 \end{center}
\end{algorithm}

\begin{proposition}[Termination, Correctness, and Completeness of \FunT{sbc}]
 Given an $n$-ary equation $c\colon f(x_1,\dots,x_n)=0$, a box $\itvK{I_1}\times\cdots\times\itvK{I_n}$ of domains, and a projection
$\itvK{c}\colon\itvK{g_i}(x)=0$ of $c$, we have:
 \begin{description}
  \item[Termination.] The call to $\FunT{sbc}(\itvK{g_i},\itvK{I_i})$ always terminates;
  \item[Correctness.] The equation $\itvK{c}$ is box consistent with respect to $x_i$ and $\FunT{sbc}(\itvK{g_i},\itvK{I_i})$;
  \item[Completeness.] No solution is lost during the tightening process:
\begin{equation*} 
\forall (r_1,\dots,r_n)\in\itvK{I_1}\times\cdots\times\itvK{I_n}\colon f(r_1,\dots,r_n)=0\implies r_i\in\FunT{sbc}(\itvK{g_i},\itvK{I_i}).
 \end{equation*}

 \end{description}
 \begin{proof} In the following, the interval \itvK{I} corresponds to the domain $\itvK{I_i}$ of $x_i$.
 
   (Termination). Algorithm~\FunT{sbc} terminates in any case since we always tighten $\itvK{I}$ in the loop 2--23 (either
   	by splitting it on Line~3, or by tightening $\itvK{I_l}$ and $\itvK{I_r}$ with Newton steps on Lines 10 and 19). Since we consider intervals
   	with floating-point bounds, of which there are finitely many, we have to reach canonicity of $\itvK{I}$ eventually. At that point, 
   		either $\itvK{g_i}(\itvK{I})$ contains $0$, and we have reached box consistency, or it does not, and 
	we can safely narrow $\itvK{I}$
   		to $\emptyset$, which both make us leave the loop.
   	
   	(Correctness). We leave the loop 2--23 if \itvK{I} is empty or if the canonical intervals at left and right bounds contain
   		solutions of \itvK{g_i} (Lines~5 and 14). In the latter case, we have the two conditions of  Eq~\eqref{eq:box-consistency} 
   	for box consistency; the former case occurs if both $\itvK{I_l}$ and $\itvK{I_r}$ do not contain a solution of $\itvK{c}$ (Lines~7 and
   		16) or if the Newton steps on Lines 10 and 19 narrow them down to $\emptyset$. By correctness of the interval Newton method,
   			this case only occurs if, once again, $\itvK{I_l}$ and $\itvK{I_r}$ do not contain a solution of $\itvK{c}$.
   	
 	(Completeness). By completeness of the Newton operator, as tightening only occurs 
either through Newton steps, or by discarding intervals that have been proved on Lines 5, 7,
14, or 16 not to contain solutions.\hfill$\Box$
 \end{proof}
\end{proposition}

For each iteration of the loop 2--23 in Algorithm~\FunT{sbc}, we have to compute
$\itvK{g}$ for intervals \itv{\lbK{\itvK{I_l}}}{\nextFloat{\lbK{\itvK{I_l}}}}, \itv{\previousFloat{\rbK{\itvK{I_r}}}}{\rbK{\itvK{I_r}}}, \itvK{I_l},
\itvK{I_r}, \itv{\lbK{\itvK{I_l}}}{\lbK{\itvK{I_l}}}, and \itv{\rbK{\itvK{I_r}}}{\rbK{\itvK{I_r}}}. We also have to compute \itvK{g'} for intervals
\itvK{I_l} and \itvK{I_r}. All these evaluations are candidates for parallelization with SIMD instructions, as presented in the next section.

\subsection{An SIMD Algorithm for Box Consistency: \FunT{vsbc}}
\label{sec:simd-algorithm}

Algorithm~\ref{alg:vsbc} presents a modification of Algorithm~\FunT{sbc} to make good use of its higher level of data parallelism thanks to the 
SIMD interval arithmetic that has been presented in Section~\ref{sec:interval-arithmetic}. 
The evaluations of \itvK{g} and \itvK{g'} are reordered to appear in pairs that can 
be evaluated in parallel. In addition, we reuse the evaluation of $\itvK{g}(\itv{\lbK{\itvK{I_l}}}{\nextFloat{\lbK{\itvK{I_l}}}})$
 and \itvK{g}(\itv{\previousFloat{\rbK{\itvK{I_r}}}}{\rbK{\itvK{I_r}}}) of Line~4 for the Newton steps of Line~22 instead of $\itvK{g}(\itv{\lbK{\itvK{I_l}}}{\lbK{\itvK{I_l}}})$ and
$\itvK{g}(\itv{\rbK{\itvK{I_r}}}{\rbK{\itvK{I_r}}})$ as was done in Algorithm~\FunT{sbc}. This choice avoids 
two evaluations of \itvK{g} at the cost of potentially slightly decreasing the tightening ability of the Newton step. The
domain computed is unaffected by this optimization. In particular, box consistency is still obtained eventually.
 	 
At each iteration of the loop between Lines 2 and 24, we perform 4 \emph{interval} evaluations of \itvK{g} and 2 \emph{interval} evaluations of \itvK{g'} for the same
cost as 2 \emph{floating-point} evaluations of $f$ and 1 \emph{floating-point} evaluation of $f'$.

\begin{algorithm}[htbp] 
 \begin{center}
  \caption{A data parallel algorithm for box consistency enforcement}
  \label{alg:vsbc}
 \algoInitCounter
 \begin{tabbing}
  \hspace*{1cm}\=12\=12\=12\=12\=12\=\kill
  \textbf{[\FunT{vsbc}]} \KWD{in:} $\itvK{g}\colon\ISet\to\ISet$; \KWD{in:} $\itvK{I}\in\ISet$\\[3pt]
   \REM{Returns the largest interval included in \itvK{I} that is box consistent w.r.t.\ \itvK{g}}\\[3pt]
 \KWD{begin}\\[3pt]
 \ACTR\>$(\Var{left\_consistent},\Var{right\_consistent})\gets(\KWD{false},\KWD{false})$\\[3pt]
 \ACTR\>\KWD{do}:\\[3pt]
 \ACTR\>\>$(\itvK{I_l},\itvK{I_r})\gets\Fun{split}(\itvK{I})$\\[3pt]
 \ACTR\>\>$(\itvK{J_l},\itvK{J_r})\gets \simd{\itvK{g}(\itv{\lbK{I_l}}{\nextFloat{\lbK{I_l}}})}{\itvK{g}(\itv{\previousFloat{\rbK{I_r}}}{\rbK{I_r}})}$\\[3pt]
 \ACTR\>\>\KWD{if} $0\notin\itvK{J_l}$: \REM{\itvK{I} not box consistent to the left?}\\[3pt]
 \ACTR\>\>\>$\itvK{I_l}\gets\itv{\nextFloat{\lbK{I_l}}}{\rbK{I_l}}$ \REM{Considering the remainder of \itvK{I_l}}\\[3pt]
 \ACTR\>\>\KWD{else}:\\[3pt]
 \ACTR\>\>\>$\Var{left\_consistent}\gets\KWD{true}$\\[3pt]
 \ACTR\>\>\KWD{if} $0\notin\itvK{J_r}$: \REM{\itvK{I} not box consistent to the right?}\\[3pt]
 \ACTR\>\>\>$\itvK{I_r}\gets\itv{\lbK{I_r}}{\previousFloat{\rbK{I_r}}}$ \REM{Considering the remainder of \itvK{I_r}}\\[3pt]
 \ACTR\>\>\KWD{else}:\\[3pt]
 \ACTR\>\>\>$\Var{right\_consistent}\gets\KWD{true}$\\[3pt]
 \ACTR\>\>\KWD{if} $\neg \Var{left\_consistent}\vee\neg\Var{right\_consistent}$:\\[3pt]
 \ACTR\>\>\> $(\itvK{K_l},\itvK{K_r}) \gets \simd{\itvK{g}(\itvK{I_l})}{\itvK{g}(\itvK{I_r})}$\\[3pt]
 \ACTR\>\>\>\KWD{if} $0\notin\itvK{K_l}:$ \REM{First checking an obvious absence of solution in \itvK{I_l}}\\[3pt]
 \ACTR\>\>\>\>$\itvK{I_l}\gets\emptyset$\\[3pt]
 \ACTR\>\>\>\KWD{if} $0\notin\itvK{K_r}:$ \REM{First checking an obvious absence of solution in \itvK{I_r}}\\[3pt]
 \ACTR\>\>\>\>$\itvK{I_r}\gets\emptyset$\\[3pt]
 \ACTR\>\>\>\REM{Performing 2 Newton steps in parallel to update both bounds}\\[1pt]
 \ACTR\>\>\>\REM{For better performances, we reuse \itvK{J_l} and \itvK{J_r}}\\
 \ACTR\>\>\>\REM{instead of $\itvK{g}(\itv{\lbK{I_l}}{\lbK{I_l}})$ and
 	 $\itvK{g}(\itv{\rbK{I_r}}{\rbK{I_r}})$}\\[3pt]
 \ACTR\>\>\>$(\itvK{I_l},\itvK{I_r})\gets\simd{\itvK{I_l}\cap\left(\itv{\lbK{I_l}}{\nextFloat{\lbK{I_l}}}-\itvK{J_l}/\itvK{g'}(\itvK{I_l})\right)}%
 	{\itvK{I_r}\cap\left(\itv{\previousFloat{\rbK{I_r}}}{\rbK{I_r}}-\itvK{J_r}/\itvK{g'}(\itvK{I_r})\right)}$\\[3pt]
 \ACTR\>\>$\itvK{I}\gets\Hull(\itvK{I_l}\cup\itvK{I_r})$ \REM{Returns the ``hull'' $\itv{\min(\lbK{I_l},\lbK{I_r})}{\max(\rbK{I_l},\rbK{I_r})}$ of the union}\\ [3pt]
 \ACTR\>\KWD{while} $\left(\itvK{I}\neq\emptyset\right)\wedge\left(\neg\Var{left\_consistent}\vee\neg\Var{right\_consistent}\right)$\\[3pt]
 \ACTR\>\KWD{return} \itvK{I}\\[3pt] 
  \KWD{end}
  \end{tabbing}
 \end{center}
\end{algorithm}

\section{Experiments}
\label{sec:experiments}

To evaluate the impact of our new algorithms, we have selected 20 instances of 12 classical test problems. Some are polynomial and others are not. The characteristics of these test problems are summarized in Table~\ref{tab:test-problems}. All problems are structurally well constrained, 
with as many equations as variables. Column ``Size'' reports the number of equations/variables. Column ``Equations'' indicates
whether all equations are polynomial (quadratic, if no polynomial has a degree greater than 2). A problem is labelled ``non-polynomial'' if at least one constraint contains a trigonometric, hyperbolic or otherwise transcendental operator. All test problems
are presented on the COPRIN web page~\cite{coprin}. 

\begin{table}[htbp]
\caption{Test problems characteristics}
\label{tab:test-problems}
 \begin{center}\small
  \begin{tabular}{l@{\hspace{10pt}}l@{\hspace{4pt}}rc}
    \hline
    Name   & Code                    & Size        & Equations \\
    \hline\hline
	Bronstein & \bench{bro} & $3$ & quadratic\\
    Broyden-banded &\bench{bb}  & $100$, $500$, and $1\,000$   & quadratic \\
    Broyden tridiagonal &\bench{bt}& $10$ and $20$          & quadratic \\
    Combustion & \bench{comb} & $10$ & polynomial\\
    Discrete Boundary Value Function &\bench{dbvf} & $10$ and $30$ & polynomial \\
    Extended Freudenstein &\bench{ef} & $30$ and $50$      & polynomial \\
    Mixed Algebraic Trigonometric & \bench{mat} & $3$ & non-polynomial\\
	Mor\'e-Cosnard &\bench{mc} & $50$ and $100$             & polynomial \\
	Robot & \bench{rob} & $8$ & quadratic \\
    Trigexp 3 &\bench{te3} & $5\,000$          & non-polynomial \\
    Troesch &\bench{tro} & $50$, $100$, and $200$                   & non-polynomial \\
    Yamamura &\bench{yam} & $6$ and $8$                    & polynomial \\
    \hline
  \end{tabular}
 \end{center}
\end{table}

All experiments were conducted on an Intel Core2 Duo T5600 1.83GHz. The Whetstone test~\cite{Curnow-Wichmann:Computing76} for 
this machine reports 1111~MIPS with a loop count equal to $100,000$. All algorithms have been implemented in an in-house C++ 
solver, with \emph{gaol}~\cite{gaol} as the underlying interval arithmetic library. The SIMD interval arithmetic presented
above has been implemented from scratch using Intel intrinsic instructions. In its current state, the library only contains 
vectorized versions of the addition, subtraction, multiplication, division, and integral power. All other SIMD
functions are emulated with sequential interval arithmetic. As a consequence, only polynomial
equation systems are entirely solved in an SIMD environment at present.  

Table~\ref{tab:experiments} reports the time spent in seconds to isolate all solutions of the test problems in domains with
a width smaller than $10^{-6}$, starting from the standard domains given on the COPRIN web page. 
An entry ``TO'' indicates a time-out (more than 30 minutes, here). Column \FunT{bc3revise} presents the results
obtained with Algorithm \FunT{bc3revise} implemented with double precision interval arithmetic on the FPU; Column \FunT{bc3vd}
corresponds to Algorithm \FunT{bc3revise} where interval arithmetic is performed in double precision with
SSE2 instructions (basic vectorization); Column \FunT{bc3vf} corresponds to Algorithm \FunT{bc3revise} where
interval arithmetic is performed in single precision with SSE2 instructions (we still perform only one interval operation
per SSE2 instruction, using only the lower half of SSE2 registers, though); Column \FunT{sbc} corresponds to Algorithm
\FunT{sbc} implemented with double precision interval arithmetic on the FPU; Column \FunT{sbcvd} corresponds to Algorithm \FunT{sbc}
where interval arithmetic is performed in double precision with SSE2 instructions; Column \FunT{vsbc} corresponds
to Algorithm \FunT{vsbc} where interval arithmetic is performed in single precision with SSE2 instructions (two interval
operations are performed in parallel); lastly, Column \FunT{vsbc\&sbcvd} corresponds to the cooperation of \FunT{vsbc} 
and \FunT{sbcvd}: \FunT{vsbc} is used until the domains are all smaller than a size fixed empirically to $0.25$; 
\FunT{sbcvd} is used afterwards.

\begin{table}[htbp]
\begin{center}
\caption{Experiments}
\label{tab:experiments}\setlength{\tabcolsep}{11pt}
\begin{tabular}{lrrrrrrr}
Problem				& \FunT{bc3revise} & \FunT{bc3vd} & \FunT{bc3vf}	&	\FunT{sbc} & \FunT{sbcvd} & \FunT{vsbc} & \FunT{vsbcvd\&sbcvd}\\	
\hline\hline
\bench{bro}			& 2.6		&	0.9 		&	0.4			&	0.3			&	0.3			&	\best{0.06}	& 0.2\\
\bench{bb} 100		& 10.4		&	3.3			&	3.2			&	2.7			&	1.2			&	\best{0.7}	& 0.8	\\
\bench{bb} 500		& 123.7		&	39.1		&	27.0		&	24.8		&	10.8		&	\best{5.3}	& 5.5	\\
\bench{bb} 1000		& 280.2		&	88.7		&	56.8		&	55.1		&	24.0		&	\best{11.1} & 11.4\\
\bench{bt} 10		& 16.9		&	5.8			&	3.4			&	2.1			&	1.0			&	\best{0.4}	& 0.5	\\
\bench{bt} 20		& 1127.4	&	382.2		&	260.3		&	141.1		&	66.7		&	\best{28.6}	& 30.4 \\
\bench{comb}		& 1.4		&	0.6			&	0.5			&	0.2			&	0.08		&	\best{0.03}	& 0.08\\
\bench{dbvf} 10		& 1.7		&	0.5			&	0.6			&	0.3			&	0.1			&	\best{0.08}	& 0.1\\
\bench{dbvf} 30		& 42.7		&	13.4		&	20.4		&	7.7			&	\best{3.2}	&	4.7		    & 3.8\nb\\
\bench{ef} 30		& 2.1		&	0.8			&	\NA			&	1.1			&	0.5			&	\NA			& \best{0.3}\\
\bench{ef} 50		& 5.1		&	1.8			&	\NA			&	2.2			&	0.9			&	\NA			& \best{0.6}\\
\bench{mat}			& 26.9		&	13.6		&	13.2		&	1.5			&	1.1			&	\best{0.4}	& 1.0\\
\bench{mc} 50		& 23.7		&	6.5			&	5.7			&	11.0		&	4.4			&	\best{3.0}	& 4.1\\
\bench{mc} 100		& 175.5		&	49.0		&	46.2		&	86.9		&	36.8		&	\best{28.0}	& 35.6\\
\bench{rob}			& 4.5		&	1.5			&	4.3			&	0.8			&	0.4			&	1.0			& \best{0.3}\\
\bench{te3} 5000	& 7.5		&	5.1			&	17.1		&	2.9			&	\best{2.7}	&	3.3			& 3.8\nb\\
\bench{tro} 50		& 24.8		&	15.3		&	29.1		&	4.4			&	3.6			&	\best{2.6}	& 3.4\\
\bench{tro} 100		& 180.1		&	112.7		&	384.2		&	30.9		&	25.4		&	33.7		& \best{24.2}\\
\bench{tro} 200		& 1341.4	&	844.4		&	\NA			&	231.1		&	188.1 		&	\NA			& \best{181}\\
\bench{yam} 6		& 14.6		&	4.3			&	4.2			&	7.3			&	2.5			&	\best{1.2}	& 1.7\\
\bench{yam} 8		& 279.1		&	84.6		&	104.0		&	91.0		&	35.1		&	\best{26.1}	& 27.7\\
\hline
\multicolumn{8}{l}{\scriptsize Times in seconds on an Intel Core2 Duo T5600 1.83GHz (whetstone $100\,000$=1111 MIPS)}\\[-3pt]
\multicolumn{8}{l}{\scriptsize Best times in bold blue. Time out \NA\ set to $1\,800\,\text{s}$.}\\[-3pt]
\multicolumn{8}{l}{\scriptsize Numbers followed by the symbol ``$\nbsymb$'' correspond to cases for which \FunT{vsbcvd\&sbc} performs 
		worse than \FunT{sbcvd} alone.}
\end{tabular}
\end{center}
\end{table}

\section{Discussion}

As can be seen from Column \FunT{sbc} of Table~\ref{tab:experiments}, enforcing box consistency
by shaving is faster than with Algorithm~\FunT{bc3revise} on all problems of our test set,
the ratio $\FunT{bc3revise}/\FunT{sbc}$ ranging from $1.9$ to $17.9$ and beyond. We also believe that \FunT{sbc} is
simpler to understand and easier to implement correctly than \FunT{bc3revise}.

Basic vectorization of interval arithmetic improves speed by up to three times (see \FunT{bc3revise} vs.\ \FunT{bc3vd}
and \FunT{sbc} vs.\ \FunT{sbcvd}) at no cost since algorithms do not have to be modified in order to benefit from it.

If we take advantage of the data parallelism inherent to Algorithm~\FunT{sbc} to vectorize interval
evaluations, leading to Algorithm~\FunT{vsbc}, we obtain even better results on all problems but \bench{dbvf 30}, 
\bench{ef 30}, \bench{ef 50}, \bench{rob}, \bench{te3 5000}, \bench{tro 100}, and \bench{tro 200}.
All other things being equal, if we emulate SIMD instructions
with double precision floating-point operations, we obtain back the times of \FunT{sbc}, which means that the  
very size of single precision floating-point numbers is the culprit here: as we vectorize 2 interval instructions with 
SSE2 registers, we must switch from double precision floats used in the rest of the program to single precision floats 
(see Figure~\ref{fig:simd-interval-addition}). The cast leads to less precision in the computation, which in turn 
has an impact on the ability to reject domains having no solutions. The same problem occurs for \FunT{bc3vf}.

\begin{figure}[h!]
\begin{center}
 \subfigure[][\bench{Extended Freudenstein} 2\label{fig:extended-freudenstein}]{\includegraphics[scale=.48]{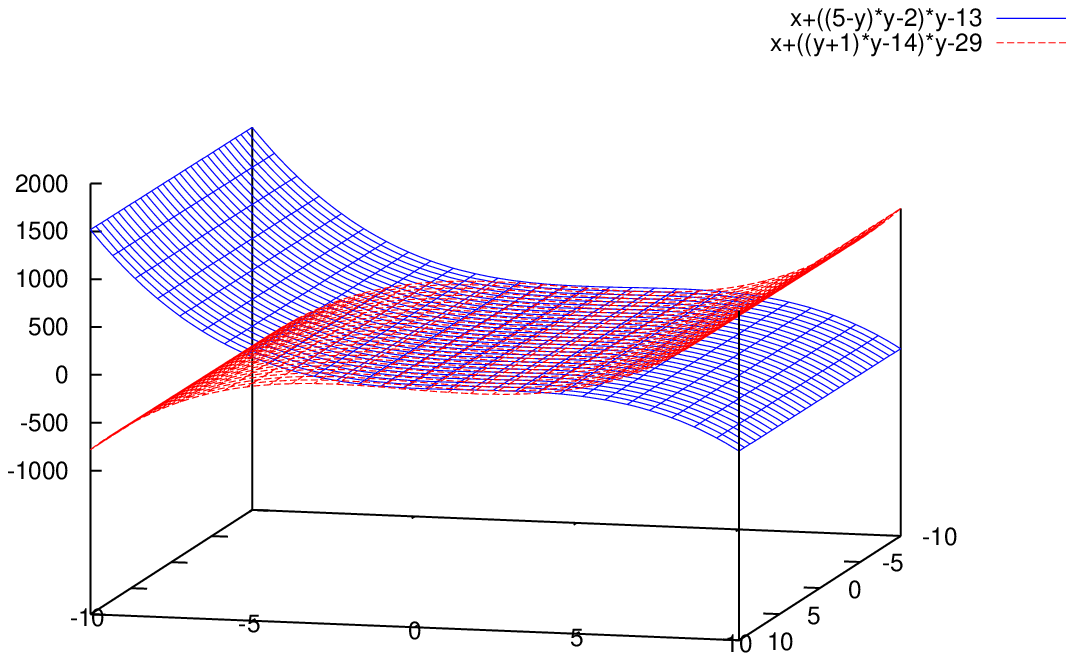}}
 \subfigure[][\bench{Troesch} 2\label{fig:troesch}]{\includegraphics[scale=.48]{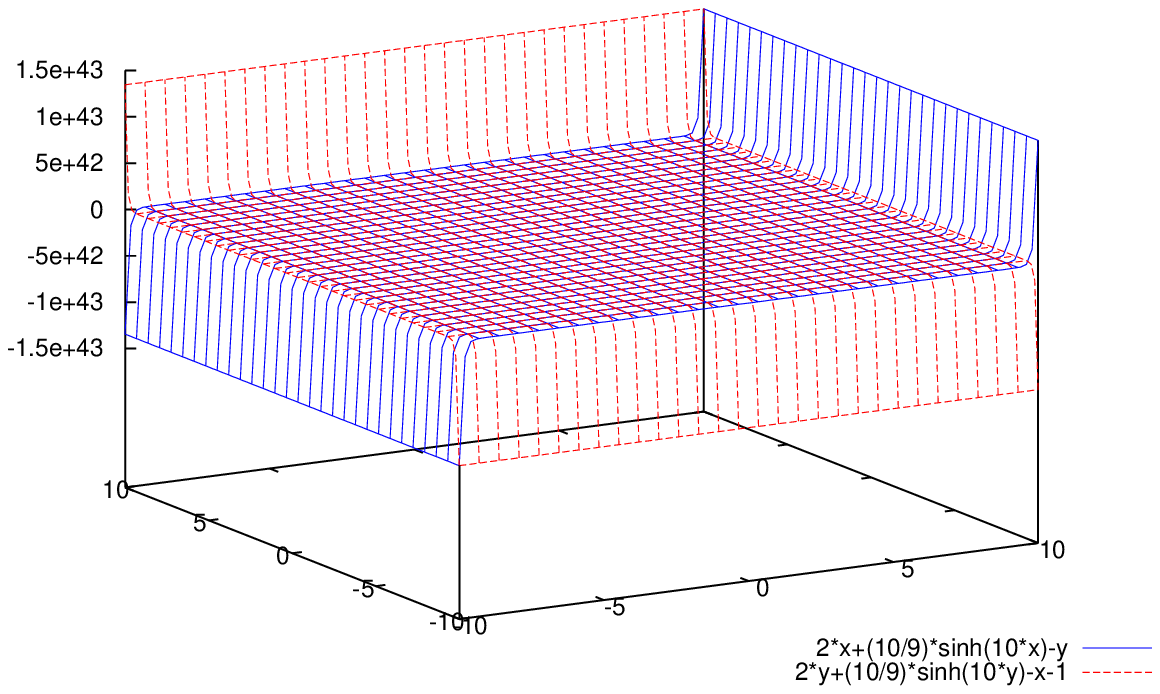}}
\caption{Slow convergence with single precision floating-point numbers}
\label{fig:ill-conditioned-problems}
\end{center}
\end{figure}

%
As a consequence, the exploration algorithm
has more branching to perform to isolate solutions. This incurs an increase in the running time that may be drastic
for ill-conditioned problems such as \bench{Troesch} or \bench{Extended Freudenstein}: as we may see in 
Figure~\ref{fig:ill-conditioned-problems} for the case of 2 equations and 2 variables, the curves for these two problems
are almost tangent to each other and to the $xy$-plane on a large surface. Each equation considered separately leads
to the computation of many quasi-zeroes that cannot be removed easily by the other equation of the problem.

There is currently no easy cure to this problem as microprocessor makers do not seem to be ready to offer SIMD instructions
on 4 double precision floats any time soon. It is still possible to quickly isolate regions of interest in ``large'' domains 
using \FunT{vsbc}, and then switch to \FunT{sbcvd} to polish the results and obtain tighter domains if necessary. 
Column \FunT{vsbc\&sbcvd} in Table~\ref{tab:experiments} shows
that this procedure indeed removes the time-out problems of \FunT{vsbc} on ill-conditioned problems, while still preserving better performances compared to \FunT{sbcvd} alone. The best cooperation scheme that maximizes performances as much as possible still remains to be found, though.

\bibliographystyle{plain}
\bibliography{3-05_Goualard-Goldsztejn}
\end{document}